\documentclass{birkjour}
 \usepackage{amsfonts,dsfont}
 \usepackage{amsthm,amsmath,amssymb}
 \usepackage[scanall]{psfrag}
 \usepackage{cleveref,enumitem}
 \newtheorem{thm}{Theorem}[section]
 \newtheorem{cor}[thm]{Corollary}
 \newtheorem{lem}[thm]{Lemma}
 
 \theoremstyle{definition}
 \newtheorem{defn}[thm]{Definition}
 \theoremstyle{remark}
 \newtheorem{rem}[thm]{Remark}
 \newtheorem{ex}{Example}
 \numberwithin{equation}{section}

 \newcommand{\R}{\mathds{R}}
 \newcommand{\N}{\mathds{N}}

 \newcommand{\m}{\mathbb{M}}
 \newcommand{\dis}[1]{\displaystyle{#1}}
 \renewcommand{\leq}{\leqslant}
 \renewcommand{\geq}{\geqslant}

\begin{document}

%
%
%
%
%
%
%
%
%

\title[Exponential stability of non-autonomous perturbed systems]
 {New results on the uniform exponential stability of non-autonomous perturbed dynamical systems}

\author[Mondher Benjemaa]{Mondher Benjemaa}

\address{%
University of Sfax\\
Laboratory Stability and Control of Systems and nonlinear PDEs\\
Soukra road, 3000, Sfax\\
Tunisia}

\email{mondher.benjemaa@fss.usf.com}

\author{Wided Gouadri}
\address{Laboratory Stability and Control of Systems and nonlinear PDEs\br
Soukra road, 3000, Sfax, Tunisia}
\email{gouadri2016@gmail.com}
\author{Mohamed Ali Hammami}
\address{Laboratory Stability and Control of Systems and nonlinear PDEs\br
Soukra road, 3000, Sfax, Tunisia}
\email{mohamedali.hammami@fss.rnu.tn}



\date{March 9, 2020}

\begin{abstract}
In this paper, we investigate the asymptotic behaviors of the solutions of nonlinear dynamic systems nearby an equilibrium point, when the nominal parts are subject to non necessarily small perturbations. We show that, under some estimates on the perturbation terms, the equilibrium point remains (globally) uniformly exponentially stable. The results we obtained can easily be applied in practice since they are based on the Gronwall-Bellman inequalities rather than the classical Lyapunov methods that require the knowledge of a Lyapunov function. Several numerical examples are presented in order to illustrate the validity of our study, especially when the standard Lyapunov approaches are useless.
\end{abstract}

\maketitle
\section{Introduction}
The study of asymptotic stability of dynamical systems is one of the most important research area in system design \cite{BH67,EL89}. In case of linear time varying (LTV) systems, it is well known that the uniform asymptotic stability of an equilibrium point is equivalent to its (global) uniform exponential stability (see \cite[Theorem 58.7]{WH67} or \cite[Theorem 4.11]{HK02} for instance). This is no longer true in case of nonlinear time varying systems. There exists in the literature several works that deal with the uniform asymptotic stability of LTV systems, see e.g. \cite{AE99,CO78,GA08,IL87,RA86,RO63,ZH16} and references therein. In concrete physical systems, perturbations (such as friction, heat, measurements, etc) occur and it is more convenient to consider perturbed LTV systems (also called the nonhomogeneous linear time varying systems). The perturbed LTV systems arise, for instance, from the linearization of a differential system around its equilibrium point, or from reducing a linear high order differential equation to a first order system, or also when using a linear feedback in the stabilization of dynamical systems, etc. In studying the effect of perturbations of these systems, it is reasonable to assume a stability property for the unperturbed or nominal system. Since the perturbations can generally be measured or at least estimated, a useful kind of stability is one for which the effect of perturbations can be studied. In this context, the Lyapunov approach gives a powerful tool to study the effects of the perturbations on the nominal system \cite{BA92,BH13,KA03,LIA93,SO01,WI67}. Indeed, the behavior of solutions under perturbations can, in some cases, be completely determined \cite{MC98,MC93,BG13,LEN98,LEN07}. The drawback of such a method is that it requires the knowledge of the Lyapunov function, which is a difficult task in general \cite{WU75}. Another approach to study the stability of perturbed systems consists of using the behavior of the solutions of the associated unperturbed system in the vicinity of the equilibrium point \cite{AB16,MM,MH15,SH10,BZ15}. Being formulated in terms of integral inequalities of Gronwall type, it is a type of stability which is easy to verify in practice, and it extends the class of systems for which the effect of perturbations can be measured. In this paper, we make use of such an approach in order to derive some conditions that ensures the (global) uniform exponential stability of perturbed systems.\\

This paper is organized as follows. In section \ref{sec_def}, we recall some definitions and integral inequalities useful to our study. Section \ref{sec_main} is dedicated to the main results of this paper, namely, we discuss the global uniform exponential stability, the uniform exponential stability and the global uniform practical exponential stability of non-autonomous perturbed systems depending on whether the perturbation term (denoted $g(t,x)$ hereafter) is linearly bounded, super-linearly bounded or sub-linearly bounded on $\|x\|$. Finally, we give in section \ref{sec_num} several numerical examples that illustrates the effectiveness of our study.

\section{Definitions and integral inequalities}\label{sec_def}
We recall in this section some definitions and integral inequalities that are useful in the sequel.
\subsection{Problem statement and definitions}
We consider the following perturbed linear time varying system (also called the nonhomogeneous linear equation \cite[\S 59]{WH67}):
\begin{equation}\label{eq4}
\dot{x}=A(t)x + g(t,x),\quad x(t_0)=x_0
\end{equation}
where $t\geq t_0\geq 0$ is the time, $x\in \R^n$ is the state, $A$ is a continuous matrix in $\m_n(\R)$ and $g : \R_+\times \R^n\longrightarrow \R^n $ is the perturbation term. We also suppose that $g(t,0)=0\ \forall\ t\geq 0$ so that the origin is an equilibrium point for the system \eqref{eq4}. Notice that equation \eqref{eq4} may arise from the linearization of the more general nonlinear system $\dot{x}=f(t,x)$ (with $A$ is the Jacobian matrix $\partial f/\partial x$ at $x=0$ in this case), or from reducing a linear high order o.d.e to a first order system, or also when using a linear feedback in the stabilization of dynamical systems, etc.
\begin{defn}
The equilibrium point $x^*=0$ is said 
\begin{itemize}
\item[(i)] uniformly exponentially stable (U.E.S) if there exist $\delta>0$, $c >0$ and $\gamma>0$ such that $\forall\ t_0\geqslant 0$, $\forall \ \Vert x_0 \Vert \leqslant \delta$,
\begin{equation*}
 \Vert x(t)\Vert \leqslant c\Vert x_0 \Vert e^{-\gamma(t-t_0)}  \quad   \forall\ t\geqslant t_0.
\end{equation*}
\item[(ii)] globally uniformly exponentially stable (G.U.E.S) if there exist $c >0$ and $\gamma>0$ such that $\forall\ t_0\geqslant 0$ and $\forall\ x_0\in \R^n$ 
\begin{equation*}
\Vert x(t) \Vert  \leqslant c\Vert x_0\Vert e^{-\gamma(t-t_0)} \quad   \forall\ t\geqslant t_0.
\end{equation*}
\end{itemize}
\end{defn}

\begin{defn}
A solution of \eqref{eq4} is said to be globally uniformly bounded if for every $\delta>0$ there exists $\eta=\eta(\delta)$ independent of $t_0$, such that $\forall\ t_0 \geqslant 0$, 
\begin{equation*}
\Vert x_0\Vert < \delta \Longrightarrow \Vert x(t)\Vert < \eta \quad \forall\ t\geqslant t_0.
\end{equation*}
\end{defn}

Let $r\geqslant 0$ and $\mathcal{B}_r=\{x\in \R^n \; /\; \Vert x\Vert \leqslant r \}$.
\begin{defn}
\begin{itemize}
\item[(i)] The ball $\mathcal{B}_r$ is uniformly stable if $\forall\ A>r$, $\exists\ \delta= \delta(A)>0$ such that $\forall\ t_0 \geqslant 0$, 
\begin{equation*}
\Vert x_0 \Vert < \delta \Longrightarrow \Vert x(t)\Vert< A \quad \forall\ t\geqslant t_0.
\end{equation*}
\item[(ii)] The ball $\mathcal{B}_r$ is globally uniformly stable if it is uniformly stable and the solutions of system \eqref{eq4} are globally uniformly bounded.
\end{itemize}
\end{defn}
\begin{defn}
\begin{itemize}
\item[(i)] The ball $\mathcal{B}_r$ is globally uniformly exponentially stable if there exist $\gamma >0$ and $c \geqslant 0$ such that $\forall \ t_0\in \R_+$ and $\forall \ x_0 \in \R^n$,
\begin{equation*}
\Vert x(t)\Vert \leqslant c \Vert x_0\Vert \,e^{-\gamma(t-t_0)}+r \quad \forall\ t\geqslant t_0.
\end{equation*}
\item[(ii)] The system \eqref{eq4} is globally practically uniformly exponentially stable if there exists $r>0$ such that $\mathcal{B}_r$ is globally uniformly exponentially stable.
\end{itemize}
\end{defn}
\subsection{Integral inequalities}
\begin{lem}(Gronwall-Bellman's inequality \cite{MH15})\label{lemma1}
Let $x$ and $\phi$ be non-negative continuous functions on $\R_+$ satisfying the inequality 
\begin{equation*}
x(t)\leq c+ \int_a^t \phi(s)x(s)\, ds ,
\end{equation*}
where $a$ and $c$ are non-negative constants. Then,
\begin{equation*}
x(t) \leq c\,e^{\int_a^t \phi(s)\,ds}\quad \forall\ t\geq a.
\end{equation*}
\end{lem}
\begin{lem}(Bihari's inequality \cite{IB56})\label{lemma2}
Let $x$ and $\phi$ be non-negative integrable functions on $[a,T]$ and let $c>0$ be an arbitrary constant. Let $\omega$ be a monotonously increasing function on $\R_+$ satisfying $\omega(0)=0$. If the inequality
\begin{align*}
x(t)\leqslant c+\int_a^t \phi(s)\, \omega(x(s))ds, \ t\in [a,T]
\end{align*}
holds, then the inequality
\begin{align*}
x(t) \leqslant \Omega^{-1}\left(\Omega(c)+\int_a^t \phi(s) ds \right), \ t\in [a,T]
\end{align*}
remains valid as long as $\Omega(c)+\dis{\int_a^t \phi(s)\, ds}$ belongs to the domain of definition of 
$\Omega^{-1}$, where the function $\Omega$ is defined by
\begin{align*}
\Omega(\eta)=\int_{\varepsilon}^{\eta}\frac{d\tau}{\omega(\tau)}, \quad \varepsilon>0, \  \eta\geqslant 0, 
\end{align*}
and $\Omega^{-1}$ is the inverse mapping of $\Omega$.
\end{lem}
\section{Main results}\label{sec_main}
The linear time varying (LTV) system
\begin{equation}\label{sys2_0}
\dot{x}=A(t)x, 
\end{equation}
is said uniformly asymptotically stable (U.A.S) if and only if there exist constants $c>0$ and $\gamma>0$ such that $\forall\ t_0\geq 0$
\begin{align}\label{RA0_ineq}
\left\| R(t,t_0)\right\| \leq c\,e^{-\gamma(t-t_0)} \quad \forall\ t\geq t_0
\end{align}
where $R(t,t_0)$ denotes the state transition matrix of the system \eqref{sys2_0} (see \cite[Theorem 4.11]{HK02} for instance). Hence, in the linear case, the uniform asymptotic stability of an equilibrium point is equivalent to its global uniform exponential stability \cite[\S 3]{CO65}.
\begin{thm}\label{theorem1}
Assume the nominal system \eqref{sys2_0} is U.A.S and suppose there exist $p\in[1,\infty)$ and a function $\phi \in L^p(\R_+,\R_+)$ such that $\|g(t,x)\|\leq \phi(t)\|x\|$ for all $(t,x)\in \R_+\times \R^n$. Then the perturbed system \eqref{eq4} is G.U.E.S. The result is still valid in case $p=\infty$ under the additional condition $\|\phi\|_\infty < \gamma/c$.
\end{thm}
\begin{proof}
Let $R(t,t_0)$ be the state transition matrix of the system \eqref{sys2_0}. It follows that the solution of \eqref{eq4} can be written as
\begin{equation}\label{x=int}
x(t) = R(t,t_0)\,x_0 + \int_{t_0}^t R(t,s)\,g(s,x(s))\,ds.
\end{equation}
Using \eqref{RA0_ineq}, we obtain
\begin{equation*}
\|x(t)\| \leq c\,e^{-\gamma(t-t_0)}\|x_0\| + c \int_{t_0}^t e^{-\gamma(t-s)}\,\phi(s)\,\|x(s)\|\,ds.
\end{equation*}
Multiply both sides by $e^{\gamma(t-t_0)}$ and denote $u(t) = \|x(t)\|\,e^{\gamma(t-t_0)}$ yields
\begin{equation}\label{ineq_u}
u(t) \leq c\,\|x_0\| + c \int_{t_0}^t \phi(s)\,u(s)\,ds.
\end{equation}
By the classical Gronwall-Bellman inequality, we deduce that
\begin{equation*}
u(t) \leq c\,\|x_0\|e^{c \int_{t_0}^t \phi(s)\,ds},
\end{equation*}
and thus
\begin{equation}\label{ineq_x}
\|x(t)\| \leq c\,\|x_0\|e^{c \int_{t_0}^t \phi(s)\,ds - \gamma(t-t_0)}.
\end{equation}
Now, if $p=1$ then 
\begin{equation*}
\|x(t)\| \leq K\,\|x_0\| e^{-\gamma(t-t_0)}
\end{equation*}
with $K = c\,e^{c \|\phi\|_1}$ and hence the system \eqref{eq4} is globally uniformly exponentially stable. Suppose $1< p < \infty$ and let $q = \frac{p}{p-1}$ be the conjugate of $p$. Using the H\"older inequality, we obtain
\begin{equation}\label{ineq_x2}
\|x(t)\| \leq c\,\|x_0\|e^{c \|\phi\|_p (t-t_0)^\frac{1}{q} - \gamma(t-t_0)}.
\end{equation}
Using that $\xi \mapsto \xi^\alpha$ is concave if $\alpha \in (0,1)$, we deduce that $\xi^\alpha \leq \varepsilon^{\alpha-1}\left(\alpha\xi + (1-\alpha)\varepsilon\right)$ for any $\varepsilon>0$. Hence, $\forall\ t\geq t_0$
\begin{equation}\label{ineq_concave}
(t-t_0)^\frac{1}{q} \leq \varepsilon^{-\frac{1}{p}}\left(\frac{1}{q}(t-t_0) + \frac{\varepsilon}{p}\right)\quad \text{ for any } \varepsilon>0.
\end{equation}
Choose $\varepsilon = \left(\frac{2c\|\phi\|_p}{q\gamma}\right)^p$ and denote $K = e^{\frac{\gamma\varepsilon}{2(p-1)}}$, we obtain from \eqref{ineq_x2} and \eqref{ineq_concave}
\begin{equation*}
\|x(t)\| \leq K\,\|x_0\|e^{-\frac{\gamma}{2}(t-t_0)}
\end{equation*}
and the system \eqref{eq4} is G.U.E.S. Finally, if $p=\infty$ then \eqref{ineq_x} implies
\begin{equation*}
\|x(t)\| \leq c\,\|x_0\|e^{-(\gamma - c \|\phi\|_\infty )(t-t_0)},
\end{equation*}
hence the system \eqref{eq4} is G.U.E.S under the condition $\|\phi\|_\infty < \gamma/c$.
\end{proof}
\begin{rem}
In the literature, the function $\phi$ is generally assumed either integrable or bounded (see e.g. \cite[Corollary 1]{ONU63} or \cite[\S 3, Theorem 6 and 9]{CO65} among others). Here, we generalize these results in case $\phi$ is $p-$integrable, with $p\in [1,\infty]$.
\end{rem}
\begin{rem}
When $p=\infty$, the condition $\|\phi\|_\infty < \frac{\gamma}{c}$ in Theorem \ref{theorem1} has also been obtained by using the Lyapunov exponent technique (see \cite[Theorem 6.1 and Theorem 6.2]{LEN08} for instance).
\end{rem}
Theorem \ref{theorem1} allows us to derive a simple and sufficient criterion on a given matrix $A(t)$ under which the linear system $\dot{x} = A(t)x$ is G.U.E.S. 
 Let us recall that the eigenvalues of $A(t)$ with negative (resp. positive) real parts for all time do not allow to conclude the exponential stability (resp. the instability) of the origin  \cite{RU96,WU74}. In view of Theorem \ref{theorem1}, we claim that if all the entries of $A(t)$, but the diagonal terms, are $p$-integrable for some $1\leq p < \infty$, and if all the diagonal terms are negatives, then the origin is G.U.E.S. More precisely, we have:
\begin{cor}\label{cor_x'=Ax}
Consider the linear time varying system \eqref{sys2_0} with $A(t) = \left(a_{ij}(t)\right)_{1\leq i,j\leq n} \in \m_n(\R)$. Suppose that
\begin{enumerate}
\item $\forall\ 1\leq i \leq n,\ \exists\ \gamma_i > 0$ and $\beta_i\geq 0$ such that $\int_{t_0}^t a_{ii}(s)\,ds \leq -\gamma_i(t-t_0) + \beta_i$, $\,\forall\ t_0\geq 0$ and $\forall\ t\geq t_0$.
\item $\exists\ p\in [1,\, \infty)$ such that $a_{ij} \in L^p(\R_+,\R)$ $\,\forall\ 1\leq i\neq j\leq n$.
\end{enumerate}
Then the origin of the system \eqref{sys2_0} is G.U.E.S.
\end{cor}
\begin{proof}
The system \eqref{sys2_0} can be written as a perturbed system
\begin{equation*}
\dot{x} = D(t)\,x + g(t,x)
\end{equation*}
with $ D = diag(a_{11},\dots,a_{nn})$ and $g(t,\cdot) = A(t) - D(t)$. It is obvious from hypothesis 1. that the origin of the nominal system $\dot{x} = D(t)\,x$ is G.U.E.S. Now, since $a_{ij} \in L^p(\R)\ \forall\ 1\leq i\neq j\leq n$ then there exists $\phi\in L^p(\R_+,\R_+)$\footnote{Since all the norms in $\m_n(\R)$ are equivalent, one can choose for instance $\phi(t) = \sqrt{n}\,\|\tilde{A}(t)\|_1$ with $\|\tilde{A}(t)\|_1 := \max_{1\leq j \leq n}\left(\sum_{\substack{i=1 \\ i\neq j}}^n |a_{ij}(t)|\right) \in L^p(\R_+,\R_+)$.} such that $\|g(t,x)\| \leq \phi(t)\|x\|$. Finally, the result follows from Theorem \ref{theorem1}.
\end{proof}
\begin{thm}\label{theorem1_1}
Assume the nominal system \eqref{sys2_0} is U.A.S and suppose there exist $p\in[1,\infty]$, a function $\phi \in L^p(\R_+,\R_+)$ and a non decreasing function $\omega\in C^1(\R_+,\R_+)$ with $\omega(0)=0$ such that $\|g(t,x)\|\leq \phi(t)\omega(\|x\|)$. Suppose either
\begin{enumerate}[leftmargin=*, label = (\alph*)]
\item $\omega'(0) = 0$
\item or $\omega'(0) > 0$ with
\begin{enumerate}[leftmargin=*, label = (\roman*)]
\item $p=1$
\item or there exists $\delta > 0$ such that $$\|\phi\|_p < \dfrac{\left(\frac{p\,\gamma}{p-1}\right)^{1-1/p}}{2^{1/p}c\|\omega'\|_{L^\infty(\mathcal{B}_\delta)}}$$ if $p\in(1,\infty]$.
\end{enumerate}
\end{enumerate}
Then the perturbed system \eqref{eq4} is uniformly exponentially stable.
\end{thm}
\begin{proof}
Let $t_0\geq 0$ and suppose the solution $x$ of \eqref{eq4} is defined on $I = [t_0,T_{\text{max}})$ with $T_{\text{max}}> t_0$. Using \eqref{RA0_ineq} and \eqref{x=int} we have for any $t\in I$
\begin{align}\label{est_xexpqomg}
\|x(t)\|e^{\gamma(t-t_0)} \leq c\,\|x_0\| + c \int_{t_0}^t e^{\gamma(s-t_0)}\phi(s)\,\omega(\|x(s)\|)\,ds.
\end{align}
Let $p\in(1,\infty]$ and let $q=\frac{p}{p-1}$ be the conjugate of $p$, with $q=1$ if $p=\infty$. Since $\omega$ is continuously derivable in $0$ with $\omega(0)=0$ then for any $R>0$ and for all $x\in \mathcal{B}_R$ we have $\omega(\|x\|) \leq L(R)\|x\|$ with $L(R)=\|\omega'\|_{L^\infty(\mathcal{B}_R)}$. Moreover, $L(R)\to 0$ as $R\to 0$ if $\omega'(0)=0$. Let $K = \max(1,2^{\frac{1}{p}}c)$ and choose $r>0$ sufficiently small such that $\gamma > \frac{1}{q}2^{q-1}\left(cL(Kr)\|\phi\|_p\right)^q$ if $\omega'(0)=0$ or $r = \delta/K$ with $\delta$ given by $(ii)$ if $\omega'(0)>0$. We will prove that the system \eqref{eq4} is uniformly exponentially stable for any initial condition $x_0\in \mathcal{B}(0,r):=\{x\in \R^n,\ \|x\|<r\}$. Suppose $\|x_0\| = r - d$ with $0<d<r$ and define the set
\begin{align*}
J = \left\{T\in I\ \ s.t.\ \ \|x(t)\|\leq Kr\ \ \forall\ t\in [t_0,T]\right\}.
\end{align*}
It is clear that $J$ is a sub-interval of $[t_0,T_{\text{max}})$ that contains $t_0$. Using \eqref{est_xexpqomg} and the H\"older inequality, we obtain for any $t\in J$
\begin{align*}
\|x(t)\|e^{\gamma(t-t_0)} & \leq c\,\|x_0\| + c \|\phi\|_p\left(\int_{t_0}^t e^{q\gamma(s-t_0)}\,\omega(\|x(s)\|)^q\,ds\right)^\frac{1}{q}\\
& \leq c\,\|x_0\| + c\,L(Kr)\,\|\phi\|_p\left(\int_{t_0}^t e^{q\gamma(s-t_0)}\,\|x(s)\|^q\,ds\right)^\frac{1}{q}.
\end{align*}
Denote $u(t) = \|x(t)\|e^{\gamma(t-t_0)}$ and $\eta = c\,L(Kr)\,\|\phi\|_p$. Using the identity 
\begin{equation}\label{id_2p}
(a+b)^\nu \leq 2^{\nu-1}(a^\nu+b^\nu)
\end{equation}
for any non negative reals $a$ and $b$ and any real $\nu \geq 1$, we obtain
\begin{align*}
u^q(t) \leq 2^{q-1}c^q\,\|x_0\|^q + 2^{q-1}\eta^q\,\int_{t_0}^t u^q(s)\,ds,
\end{align*}
yielding by the classical Gronwall-Bellman inequality
\begin{align*}
u^q(t) \leq 2^{q-1}c^q\,\|x_0\|^q\,e^{2^{q-1}\eta^q\,(t-t_0)},
\end{align*}
or equivalently
\begin{align}\label{est_xexpq_finomg}
\|x(t)\| & \leq 2^{\frac{1}{p}}c\,\|x_0\|\,e^{-\left(\gamma-\frac{1}{q}2^{q-1}\eta^q\right)(t-t_0)}\nonumber \\
& \leq K\|x_0\|e^{-\left(\gamma-\frac{1}{q}2^{q-1}\eta^q\right)(t-t_0)}\quad \forall\ t\in J.
\end{align}
In case $p=1$, then one can choose $r = \delta/K$ where $\delta>0$ is arbitrary and $K=\max(1,ce^{cL(\delta)\|\phi\|_1})$. Using similar arguments as previously, one can show that 
\begin{align}\label{est_xexpq_finomg_p=1}
\|x(t)\| & \leq c\,\|x_0\|\,e^{cL(Kr)\|\phi\|_1}e^{-\gamma(t-t_0)} \nonumber \\
& \leq K\|x_0\|e^{-\gamma(t-t_0)} \quad \forall\ t\in J.
\end{align}
Now we prove that $J=I$. Suppose by contradiction that $J\subsetneq I$ then $J = [t_0,T_1)$ with 
$$T_1 : = \sup\left\{T\in I\ \ s.t.\ \ \|x(t)\|\leq K r\ \ \forall\ t\in [t_0,T]\right\} < T_{\text{max}}.$$ Tacking the limit as $t$ tends to $T_1$ in \eqref{est_xexpq_finomg} or \eqref{est_xexpq_finomg_p=1} and using the continuity of $x$, one deduce that $\|x(T_1)\| \leq K\|x_0\| = K(r-d)$ and $T_1\in J$. Using again the continuity of $x$ in $T_1$, one can find $0 < \nu < T_{\text{max}} - T_1$ such that for any $t\in [T_1,T_1+\nu]$ we have
\begin{align*}
\|x(t)\| \leq \|x(T_1)\| + Kd \leq Kr
\end{align*}
which is in contradiction with the definition of $T_1$. Hence $T_1 = T_{\text{max}}$ and $J=I$. Finally, using \eqref{est_xexpq_finomg} one deduce that $T_{\text{max}} = +\infty$ and the origin of the system \eqref{eq4} is uniformly exponentially stable.
\end{proof}
\begin{cor}\label{cor_1}
Assume the nominal system \eqref{sys2_0} is U.A.S and suppose there exist $p\in[1,\infty]$ and a function $\phi \in L^p(\R_+,\R_+)$ such that $\|g(t,x)\|\leq \phi(t)\,\|x\|^\ell$ with $\ell>1$. Then the system \eqref{eq4} is uniformly exponentially stable.
\end{cor}
\begin{proof}
A direct consequence of Theorem \ref{theorem1_1}.
\end{proof}
\begin{rem}
If $\omega(0) = \omega'(0) = 0$ then $\omega(x) = o(x)$ and the statement $(a)$ in Theorem \ref{theorem1_1} is equivalent to the Lyapunov theorem for quasi-linear system (see e.g. \cite[Theorem 56.2]{WH67}).
\end{rem}
\begin{rem}
The results of Theorem \ref{theorem1_1} are not guaranteed if the asymptotic stability of the nominal system \eqref{sys2_0} is not uniform. A counter example can be found in the Perron effects \cite{PE30}.
\end{rem}
\begin{rem}\label{remark_notexp}
Generally, one can not expect a global exponential stability in Theorem \ref{theorem1_1}. Consider for instance the non linear scalar equation $\dot{x}=-x+2e^{-t}x^2$. It is not difficult to show that the solution of this equation with the initial condition $x(0)=x_0\in \R$ is given by $$x(t) = \dfrac{x_0}{x_0e^{-t}-(x_0-1)e^t}.$$
Hence, for any $\mu\geq 1$, if $x_0\leq 1-1/\mu$ then $|x(t)| < \mu |x_0|e^{-t}$ $\forall\ t\geq 0$ and the solution is exponentially stable. However, if $x_0= 1$ then $x(t) = e^t$ and the solution is not bounded.
\end{rem}
The following theorem gives a result on the practical stability of the system \eqref{eq4}. Thus, the hypothesis $g(t,0) = 0$ for all time $t$ is not required in what follows.
\begin{thm}\label{theorem2_-1}
Assume the nominal system \eqref{sys2_0} is U.A.S and suppose that $$\|g(t,x)\|\leq \phi(t)\,\omega(\|x\|) + \lambda(t)$$ where
\begin{itemize}
\item $\phi \in L^p(\R_+,\R_+)$ with $p\in[1,\infty]$.
\item $\omega : \R_+\rightarrow \R_+$ is a non decreasing continuous function s.t. $\lim_{x\to \infty}\frac{\omega(x)}{x^\ell}=0$ with $\ell\in (0,1]$ if $p=1$ or $\ell\in (0,1)$ if $p\in(1,\infty]$.
\item $\lambda \in L^r(\R_+,\R_+)$ with $r\in [1,\infty]$ or $\lim_{t\to \infty}\lambda(t)=0$.
\end{itemize}
Then
\begin{itemize}
\item[(i)] The solution $x$ of \eqref{eq4}, if there exists, is defined on $[t_0,\infty)$.
\item[(ii)] The system \eqref{eq4} is globally practically uniformly exponentially stable.
\item[(iii)] If $p$ and $r\in [1,\infty)$ then $\lim_{t\to \infty}\|x(t)\|=0$.
\end{itemize}
\end{thm}
\begin{proof}
Let us first notice that it is proven in \cite{ST67} that if
\begin{equation}\label{cond_t_t+1}
\lim_{t\to \infty}\int_t^{t+1}\lambda(s)\,ds = 0
\end{equation} then 
\begin{equation}\label{cond_exp-t}
\lim_{t\to \infty}e^{-\gamma t}\int_0^t e^{\gamma s}\lambda(s)\,ds = 0
\end{equation}
for any $\gamma > 0$ (see \cite[Lemmas 3.4, 3.5 and 3.6]{ST67}). Obviously, if $\lambda \in L^r(\R_+,\R_+)$ with $r\in [1,\infty)$ or if $\lim_{t\to \infty}\lambda(t)=0$ then the condition \eqref{cond_t_t+1} is satisfied, and so is \eqref{cond_exp-t}. Hence there exists $M\geq 0$ such that $\int_{0}^t e^{-\gamma(t-s)}\lambda(s)\,ds \leq M$ for all $t\geq 0$ and for any $r\in[1,\infty]$ such that $\lambda \in L^r(\R_+,\R_+)$. On another hand, since $\lim_{x\to \infty}\omega(x)/x^\ell=0$ and $\omega$ is non decreasing then there exists $A\geq 0$ sufficiently large such that $\omega(x) \leq x^\ell + \omega(A)$ for all $x\geq 0$. Using \eqref{RA0_ineq} and \eqref{x=int} we obtain for any $t\geq t_0$
\begin{align}\label{ineqthm2}
\|x(t)\|e^{\gamma(t-t_0)} & \leq c\,\|x_0\| + c \int_{t_0}^t e^{\gamma(s-t_0)}\,\lambda(s)\,ds + c \int_{t_0}^t e^{\gamma(s-t_0)}\,\phi(s)\,\omega(\|x(s)\|)\,ds \nonumber\\
& \leq c\,\|x_0\| + c \int_{t_0}^t e^{\gamma(s-t_0)}\big{(}\lambda(s)+\omega(A)\,\phi(s)\big{)}ds \nonumber \\
& \quad + c \int_{t_0}^t e^{\gamma(s-t_0)}\,\phi(s)\|x(s)\|^\ell\,ds,
\end{align}
or equivalently
\begin{equation}\label{ineqthm2_u}
u(t) \leq h(t) + c \int_{t_0}^t \phi(s)\,e^{(1-\ell)\gamma(s-t_0)}\,u(s)^\ell\,ds
\end{equation}
where $u(t) = \|x(t)\|e^{\gamma(t-t_0)}$ and $h(t)=c\,\|x_0\| +c \int_{t_0}^t e^{\gamma(s-t_0)}\big{(}\lambda(s)+\omega(A)\,\phi(s)\big{)}\,ds$. Now we discuss two cases.\\

\noindent $\boldsymbol{1}^\text{\bf st}$ \textbf{case} : $p\in(1,\infty]$ and $\ell \in (0,1)$. It follows using the Bihari inequality and the identity \eqref{id_2p} that
\begin{align}\label{u(t)maj1}
u(t) & \leq \left(h(t)^{1-\ell} + c(1-\ell)\int_{t_0}^t \phi(s)\,e^{(1-\ell)\gamma(s-t_0)}\,ds\right)^\frac{1}{1-\ell} \nonumber \\
& \leq 2^\frac{\ell}{1-\ell}\left(h(t) + \left(c(1-\ell)\int_{t_0}^t \phi(s)\,e^{(1-\ell)\gamma(s-t_0)}\,ds\right)^\frac{1}{1-\ell}\right).
\end{align}
Using the H\"older inequality, we obtain
\begin{align}\label{u(t)maj2}
u(t)& \leq 2^\frac{\ell}{1-\ell}\left(h(t) + \left(\dfrac{c^q(1-\ell)^q\|\phi\|_p^q}{(1-\ell)\gamma q}\right)^\frac{1}{(1-\ell)q}\left(e^{(1-\ell)\gamma q(t-t_0)}-1\right)^\frac{1}{(1-\ell)q}\right) \nonumber \\
& \leq 2^\frac{\ell}{1-\ell}\left(h(t) + \left(c\|\phi\|_p\gamma^{-\frac{1}{q}}\right)^\frac{1}{1-\ell}e^{\gamma(t-t_0)}\right)
\end{align}
where $q=\frac{p}{p-1}$ is the conjugate of $p$ (with $q=1$ if $p=\infty$).
It follows that $\forall\ t\geq t_0$
\begin{align*}
\|x(t)\| & \leq 2^\frac{\ell}{1-\ell}c\left(\|x_0\|e^{-\gamma(t-t_0)} + \int_{t_0}^t e^{-\gamma(t-s)}\big{(}\lambda(s)+\omega(A)\phi(s)\big{)}\,ds\right) \\
& \quad + \left(2^\ell c\|\phi\|_p\gamma^\frac{1-p}{p}\right)^\frac{1}{1-\ell} \\
& \leq K\|x_0\|e^{-\gamma(t-t_0)} + \delta
\end{align*}
with $K = 2^\frac{\ell}{1-\ell}c$ and $\delta = K(1+\omega(A))M + \left(2^\ell c\|\phi\|_p\gamma^\frac{1-p}{p}\right)^\frac{1}{1-\ell}$, which means that the ball $\mathcal{B}_\delta$ is globally uniformly exponentially stable and the solution $x$, since bounded, is defined on $[t_0,+\infty)$.\\

\noindent $\boldsymbol{2}^\text{\bf nd}$ \textbf{case} : $p=1$ and $\ell \in (0,1]$. In fact, we only need to treat the case $\ell = 1$ since if $\ell \in (0,1]$ then $x^\ell \leq x$ for all $x\geq 1$ and hence $\omega(x) \leq x + \omega(A)$ if we choose $A\geq 1$ (recall $\omega$ is non decreasing). Using the classical Gronwall-Bellman inequality, we obtain from \eqref{ineqthm2}
\begin{equation*}
u(t) \leq h(t)\,e^{c\int_{t_0}^t \phi(s)\,ds} \leq h(t)\,e^{c\|\phi\|_1}
\end{equation*}
where $u(t) = \|x(t)\|e^{\gamma(t-t_0)}$ and $h(t)=c\,\|x_0\| +c \int_{t_0}^t e^{\gamma(s-t_0)}\big{(}\lambda(s)+\omega(A)\,\phi(s)\big{)}\,ds$. Multiplying both sides by $e^{-\gamma(t-t_0)}$ we obtain
\begin{align*}
\|x(t)\| & \leq \left(c\,\|x_0\|e^{-\gamma(t-t_0)} +c \int_{t_0}^t e^{-\gamma(t-s)}\big{(}\lambda(s)+\omega(A)\,\phi(s)\big{)}\,ds\right)\,e^{c\|\phi\|_1} \\
& \leq K\|x_0\|e^{-\gamma(t-t_0)} + \delta
\end{align*}
with $K = c\,\,e^{c\|\phi\|_1}$ and $\delta = K(1+\omega(A))M$, which means that system \eqref{eq4} is globally practically uniformly exponentially stable and the solution $x$ is defined on $[t_0,+\infty)$.\\

Now we prove that the origin is attractive. Let $R>0$ such that $\Vert x(t)\Vert \leqslant R \ \ \forall\ t\geqslant t_0$ and denote $\varphi(t) = c \Vert x_0\Vert e^{-\gamma(t-t_0)} + c\int_{t_0}^t e^{-\gamma(t-s)}\big{(}\lambda(s)+\omega(A)\phi(s)\big{)}\,ds$. It follows using \eqref{ineqthm2}
\begin{align*}
\Vert x(t)\Vert &\leqslant \varphi(t) +c \int_{t_0}^t e^{-\gamma(t-s)}\,  \phi(s)\,  R^\ell\, ds\\
& \leq \varphi(t) + c\,R^\ell \left(\int_{t_0}^t e^{-\gamma(t-s)}\,\phi^p(s)\, ds\right)^\frac{1}{p}\left(\int_{t_0}^t e^{-\gamma(t-s)}\,ds\right)^\frac{1}{q}\\
& \leq \varphi(t) + c\gamma^{-\frac{1}{q}}\,R^\ell\left(\int_{t_0}^t e^{-\gamma(t-s)}\,  \phi^p(s)\, ds\right)^\frac{1}{p}.
\end{align*}
with $q=\frac{p}{p-1}$ is the conjugate of $p$ ($q=\infty$ if $p=1$). Using the identity \eqref{id_2p} 
, we obtain
\begin{align*}
2^{1-p}\Vert x(t)\Vert^p & \leq \varphi^p(t)+ K \int_{t_0}^t e^{-\gamma(t-s)}\,  \phi^p(s)\, ds \\
& \leqslant \varphi^p(t) + K \int_0^{t/2} e^{-\gamma(t-s)}\,\phi^p(s)\,ds + K \int_{t/2}^t e^{-\gamma(t-s)}\, \phi^p(s)\, ds \\
&\leqslant \varphi^p(t) + K \,e^{-\gamma t/2} \int_0^{t/2} \phi^p(s)\, ds + K  \int_{t/2}^t \phi^p(s)\, ds
\end{align*}
with $K = c^p\gamma^{1-p}R^{\ell p}$ is a positive constant. Finally, the result follows directly by noticing that $\int_0^{\infty} \phi^p(s)\,ds$ is finite and $\lim_{t\to \infty} \varphi(t) = 0$ if $p$ and $r$ are in $[1,\infty)$.
\end{proof}

\begin{cor}\label{cor_2}
Assume the nominal system \eqref{sys2_0} is U.A.S and suppose there exist $p\in[1,\infty]$ and a function $\phi \in L^p(\R_+,\R_+)$ such that $\|g(t,x)\|\leq \phi(t)\,\|x\|^\ell+\lambda(t)$ with $\ell\in(0,1)$ and $\lim_{t\to \infty}\lambda(t)=0$ or $\lambda \in L^r(\R_+,\R_+)$ with $r\in [1,\infty]$. Then the system \eqref{eq4} is globally practically uniformly exponentially stable. Moreover, if $p<\infty$ and $r<\infty$ then $\lim_{t\to \infty}\|x(t)\|=0$.
\end{cor}
\begin{proof}
A direct consequence of Theorem \ref{theorem2_-1}.
\end{proof}

%
%
\section{Numerical results}\label{sec_num}
In what follows, we give some numerical examples to illustrate our theoretical study. The first example is relative to Theorem \ref{theorem1}. The second and third examples illustrate Theorem \ref{theorem1_1} and the fourth example is relative to Theorem \ref{theorem2_-1}. All the illustrations have been performed with the software Matlab.
%
%
\begin{ex}\label{example2}
Consider the linear system
\begin{align} \label{sys_example2}
\begin{cases}
\dot{x}_1=(-1+2\cos(t))x_1+\frac{x_2}{\sqrt{1+t^3}}\\
\dot{x}_2=(-1-2\cos(t))x_2+\frac{x_1}{\sqrt{1+t^3}}
\end{cases}
\end{align}
with $x=(x_1,x_2)\in \R^2$ is the state of the system. The system \eqref{sys_example2} can be written as
\begin{equation*}
\dot{x}=A(t)x
\end{equation*}
with $$A(t) = \begin{pmatrix}
-1+2\cos(t) & \frac{1}{\sqrt{1+t^3}} \\ \frac{1}{\sqrt{1+t^3}} & -1-2\cos(t)
\end{pmatrix}.$$ Since $\forall\ t\geq t_0 \geq 0$
$$\int_{t_0}^t a_{ii}(s)\,ds \leq -(t-t_0) + 2\left|\sin(t) - \sin(t_0)\right| \leq -(t-t_0) + 4,\quad i=1,2$$ and $a_{12}(t) = a_{21}(t) = \frac{1}{\sqrt{1+t^3}} \in L^1(\R_+,\R_+)$, then all the hypothesis of Corollary \ref{cor_x'=Ax} are fulfilled and one conclude that the system \eqref{sys_example2} is G.U.E.S. The figure \ref{fig0} shows the time evolution of the states $(x_1,\,x_2)$ and $\log\|x\|$ of the system \eqref{sys_example2} with the initial states $(x_{1,0},\,x_{2,0}) = (3,1)$. One can notice that the mapping $t \mapsto \log\|x(t)\|$ remains below a line with negative slope for all time $t$, which confirms that the origin is exponentially stable as predicted by theory. 
\begin{figure}
\begin{center}
\includegraphics[scale=0.3]{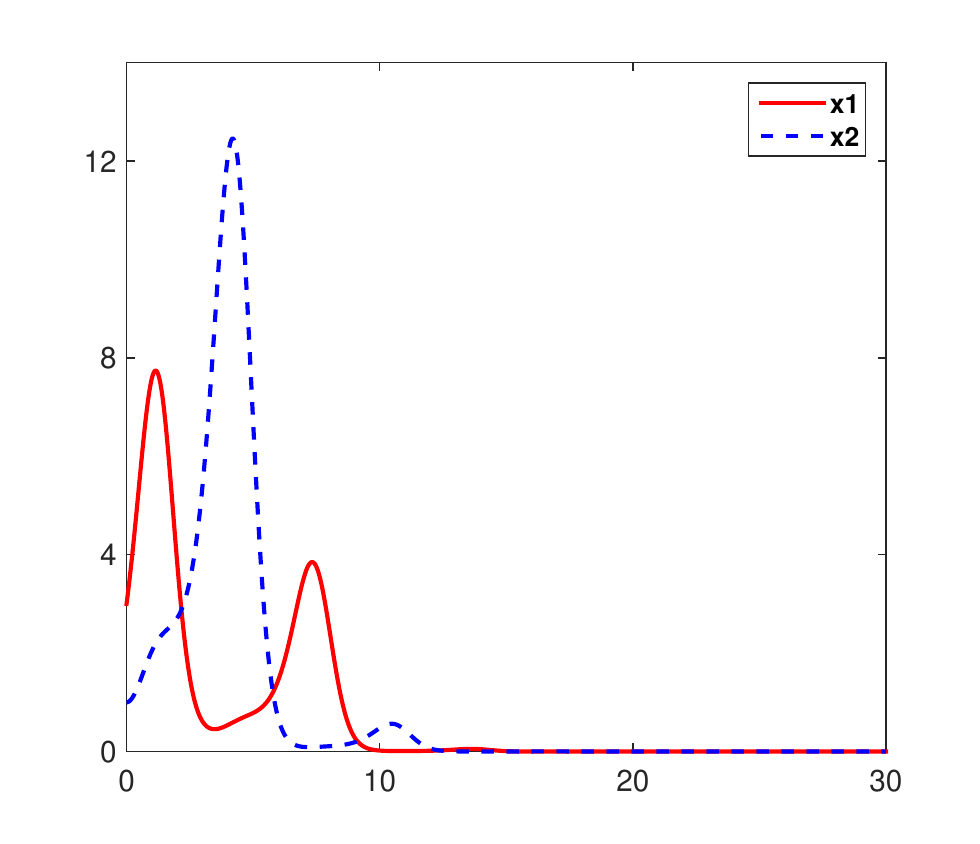} \includegraphics[scale=0.3]{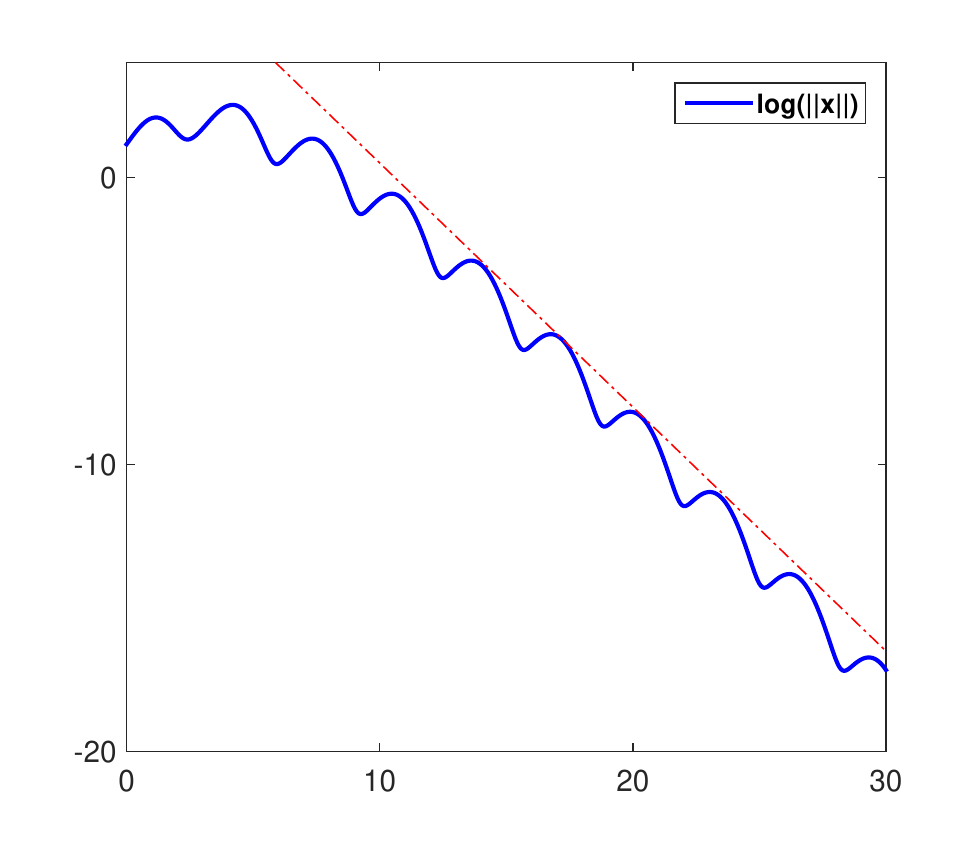}
\caption{Time evolution of the states $(x_1,x_2)$ (left) and $\log \|x\|$ (right) of the system \eqref{sys_example2}. The parameter used are $t_0=0$ and $x_0=(3,1)$. The negative slope of the dashed line (right figure) indicates the exponential decay of $\|x(t)\|$ with respect to the time $t$.\label{fig0}}
\end{center}
\end{figure}
Let us also remark that the Lyapunov function relative to the time invariant system fails to be an adequate candidate function to the system \eqref{sys_example2}. Consequently, one cannot conclude the global uniform exponential stability of the origin by using such a function. Indeed, let $V(x) = \frac{1}{2}\|x\|^2$ be the Lyapunov function of the time invariant system $\dot{x}=-x$, then the derivative of $V$ along the solutions of \eqref{sys_example2} is given by
\begin{align*}
\dot{V}(x) = \langle x,\, \dot{x}\rangle = -x_1^2-x_2^2+2\cos(t)\left(x_1^2-x_2^2\right)+\frac{2x_1 x_2}{\sqrt{1+t^3}}.
\end{align*}
It follows that if $x\in D=\{(x_1,x_2)\in \R^2,\ x_1=x_2\}$ then
\begin{equation*}
\dot{V}(x) = \left(\frac{2}{\sqrt{1+t^3}}-1\right)\|x\|^2
\end{equation*}
which is positive for $0\leq t < \sqrt[3]{3}$ and $\dot{V}$ is not negative definite for all $t\geq 0$.
\end{ex}
%
%
\begin{ex}Consider the following non linear system
\begin{align}\label{system_sin} 
\begin{cases}
\dot{x}_1=-x_1-t\,x_2 + \phi(t)\,x_2\,\sqrt[3]{x_1}\,e^{\|x\|}\\
\dot{x}_2=-x_2+t\,x_1 + \phi(t)\,x_1\,\sqrt[3]{x_2}\, e^{\|x\|}.
\end{cases}
\end{align}
with $\phi\in L^1(\R_+)$ is given over each interval $[n,\, n+1]$, $n\in \N$ by 
\begin{align*}
\phi_{|_{[n,\,n+1]}}(t) = \left\{
\begin{array}{ll}
2n^3(n-h(t_{1,n}))\left(t-n-\frac{1}{2}\right)+n & \text{ if }\ t\in\left[t_{1,n},\, n+\frac{1}{2}\right] \\ & \\
-2n^3(n-h(t_{2,n}))\left(t-n-\frac{1}{2}\right)+n & \text{ if }\ t\in\left[n+\frac{1}{2},\, t_{2,n}\right]\\ & \\
h(t) & \text{ if }\ t\in [n,\, t_{1,n}] \cup [t_{2,n},\, n+1]
\end{array}
\right.
\end{align*}
where $t_{1,n} = n+\frac{1}{2}-\frac{1}{2n^3}$, $t_{2,n} = n+\frac{1}{2}+\frac{1}{2n^3}$ and $h(t) = \frac{1}{1+t^2}$ (see figure \ref{fig_varphi}).
\begin{figure}
\begin{center}
\includegraphics[scale=0.4]{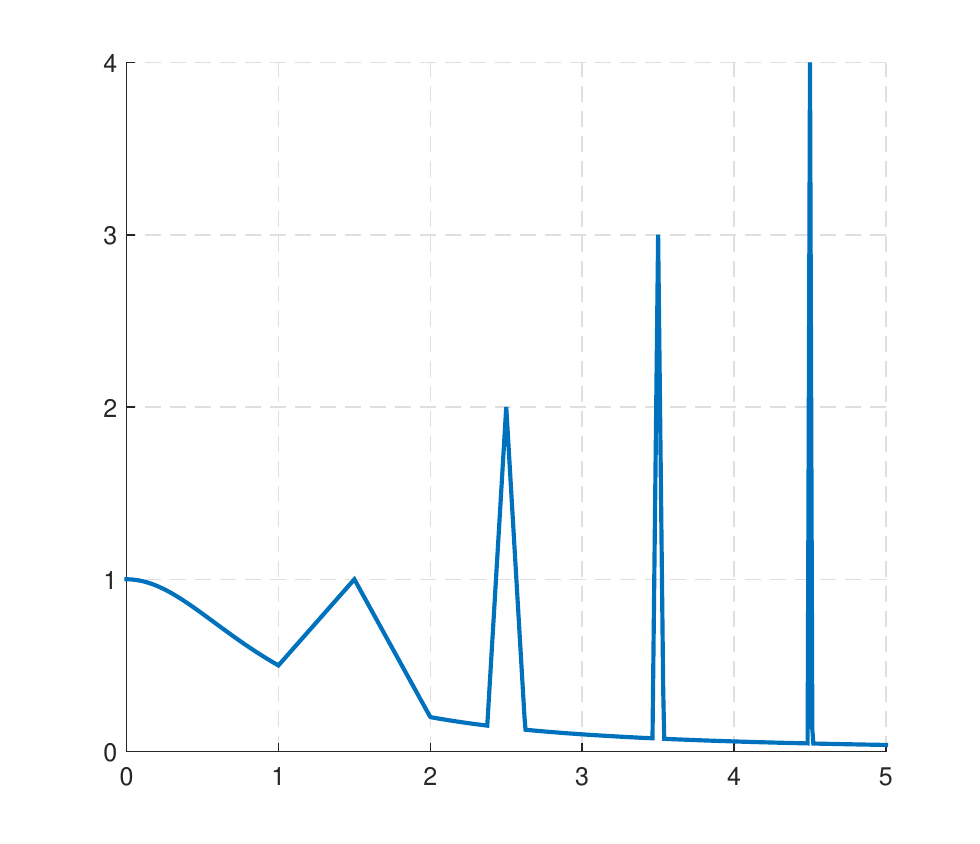}
\caption{$\phi(t)$. \label{fig_varphi}}
\end{center}
\end{figure}
The system \eqref{system_sin} can be written as
\begin{equation*}
\dot{x}=A(t)x+g(t,x),
\end{equation*}
where $x=(x_1,x_2)\in \mathbb{R}^2$ is the state of the system, 
\begin{equation*}
A(t)=\begin{pmatrix} -1 & -t \\ t & -1 
\end{pmatrix}\quad and \quad g(t,x) = \phi(t)\big{(}x_2\,\sqrt[3]{x_1},x_1\,\sqrt[3]{x_2}\big{)}e^{\|x\|}.
\end{equation*}
The transition matrix of the nominal system is given by
\begin{equation*}
R_A(t,t_0) = e^{-(t-t_0)}\begin{pmatrix}
\cos\left(\zeta(t,t_0)\right) & -\sin\left(\zeta(t,t_0)\right) \\
\sin\left(\zeta(t,t_0)\right) & \cos\left(\zeta(t,t_0)\right) 
\end{pmatrix}
\end{equation*}
with $\zeta(t,t_0) = \frac{1}{2}(t^2-t_0^2)$. It follows that $\Vert R_A(t,t_0)\Vert =e^{-(t-t_0)}$ and the unperturbed system is globally uniformly exponentially stable. Using the Young identity $ab\leq \frac{1}{4}\left(a^4 + 3b^{4/3}\right)$ for any $(a,b)\in \R_+^2$ and the inequality $1+z^\alpha\leq (1+z)^\alpha $ for any $z\geq 0$ and $\alpha\geq 1$, we obtain
\begin{align*}
\|g(t,x)\|^2 & = \phi^2(t)\left(x_1^2\,x_2^{2/3}+x_2^2\,x_1^{2/3}\right)e^{2\|x\|}\\
& \leq \phi^2(t)\left(x_1^{8/3}+x_2^{8/3}\right)e^{2\|x\|}\\
& \leq \phi^2(t)\|x\|^{8/3}e^{2\|x\|}
\end{align*}
and thus $\|g(t,x)\|\leq \phi(t)\,\omega(\|x\|)$ with $\omega(x) = x^{4/3}e^x$. Theorem \ref{theorem1_1} stipulates that the system \eqref{system_sin} is uniformly exponentially stable. The figure \ref{fig2} shows the time evolution of the states $(x_1,\,x_2)$ and $\log\|x\|$ of the system \eqref{system_sin} with the initial states $(x_{1,0},\,x_{2,0}) = (0.3,0.6)$. One can notice that the origin is exponentially stable as predicted by theory. Let us remark that, as in the previous example, the Lyapunov function relative to the nominal system fails to be an adequate candidate function to the perturbed system \eqref{system_sin}. Indeed, the derivative of $V(t,x) := \frac{1}{2}\|x\|^2$ along the trajectories of \eqref{system_sin} is given by
\begin{align*}
\dot{V}(x) = -x_1^2-x_2^2+\phi(t)\,x_1x_2\left(\sqrt[3]{x_1}+\sqrt[3]{x_2}\right)e^{\|x\|}.
\end{align*}
It follows that if $x\in D=\{(x_1,x_2)\in \R^2,\ x_1=x_2\}$ then
\begin{equation*}
\dot{V}(x) = -2x_1^2\left(1-\phi(t)\sqrt[3]{x_1}e^{x_1\sqrt{2}}\right).
\end{equation*}
Since $\phi$ is unbounded over $\R_+$ then $\dot{V}$ can not be negative definite for all $t\geq 0$ even for small enough $x\in D$.
\begin{figure}
\begin{center}
\includegraphics[scale=0.3]{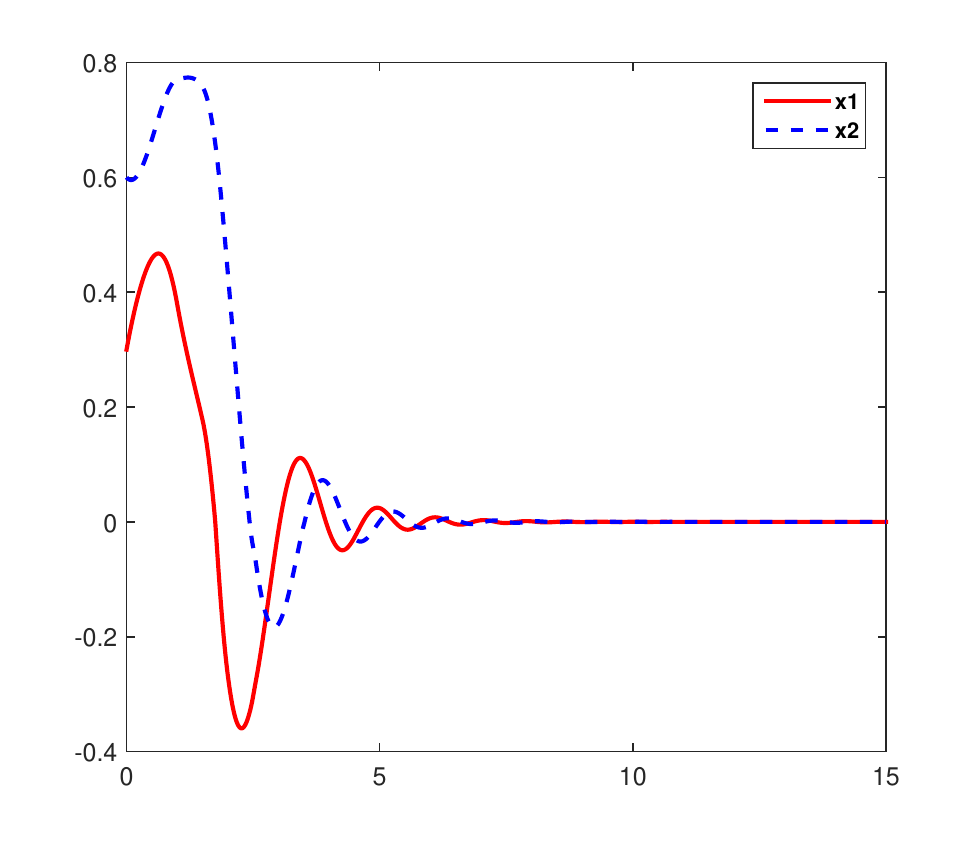} \includegraphics[scale=0.3]{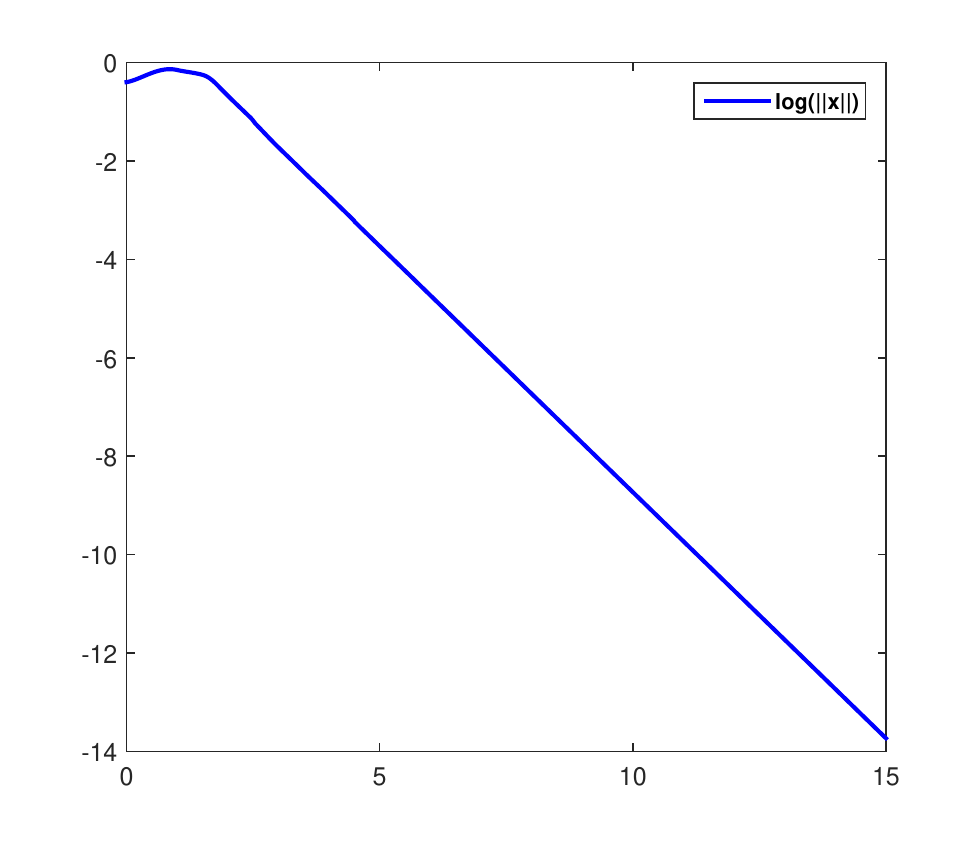}
\caption{Time evolution of the states $(x_1,x_2)$ (left) and $\log \|x\|$ (right) of the system \eqref{system_sin}. The parameter used are $t_0=0$ and $x_0=(0.3,0.6)$. \label{fig2}}
\end{center}
\end{figure}
\end{ex}
%
%
\begin{ex}
Let $\alpha \neq 0$ and consider the following non linear system
\begin{align}\label{system_exp-1} 
\begin{cases}
\dot{x}_1=-x_1 + \alpha\left(e^{\|x\|}-1\right)\\
\dot{x}_2=-x_2 - \alpha\left(e^{\|x\|}-1\right)
\end{cases}
\end{align}
which can be written as
\begin{equation*}
\dot{x}=Ax+g(t,x),
\end{equation*}
where $x=(x_1,x_2)\in \mathbb{R}^2$ is the state of the system, $A=-I_2$ and $g(t,x) = \alpha\left(e^{\|x\|}-1\right)(1,-1)$. The unperturbed nominal system is clearly globally uniformly exponentially stable. Moreover, we have $$\|g(t,x)\| = |\alpha|\sqrt{2}\left(e^{\|x\|}-1\right) = \phi(t)\omega(\|x\|)$$
with $\phi(t)=|\alpha|\sqrt{2}$ and $\omega(x)=e^x-1$. Since $\omega'(0) = 1>0$ then Theorem \ref{theorem1_1} asserts that the system \eqref{system_exp-1} is uniformly exponentially stable provided that there exist $\delta > 0$ such that $|\alpha|\sqrt{2} < e^{-\delta}$. It follows that for any $|\alpha| < 1/\sqrt{2}$ the system \eqref{system_exp-1} is U.E.S. The figure \ref{fig3} shows the time evolution of the states $(x_1,\,x_2)$ and $\log\|x\|$ of the system \eqref{system_exp-1} with the initial states $(x_{1,0},\,x_{2,0}) = (0.5,1)$ for $\alpha = 0.5$. One can notice that the solution decays exponentially as the time increases. On another hand, the figure \ref{fig4} shows the time evolution of the states $(x_1,\,x_2)$ of the system \eqref{system_exp-1} with the initial states $(x_{1,0},\,x_{2,0}) = (10^{-10},10^{-9})$ when $\alpha = 1$. One can notice that the solutions are not bounded and thus the origin of the system \eqref{system_exp-1} is not exponentially stable. This is in total agreement with Theorem \ref{theorem1_1} and Remark \ref{remark_notexp}.
\begin{figure}
\begin{center}
\includegraphics[scale=0.3]{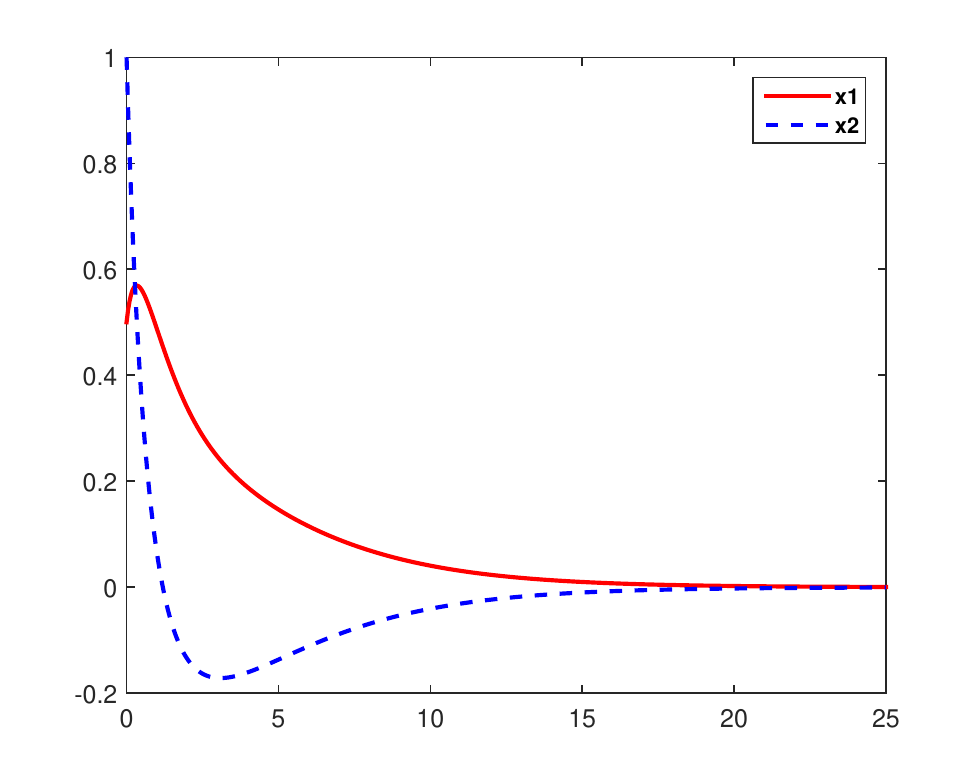} \includegraphics[scale=0.3]{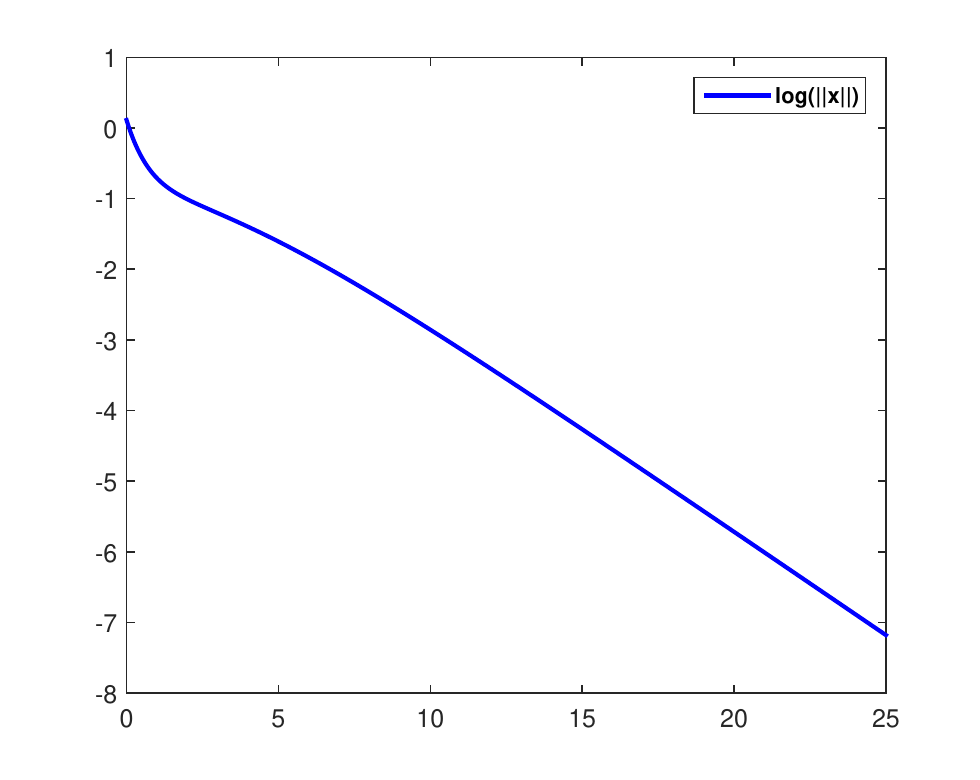}
\caption{Time evolution of the states $(x_1,x_2)$ (left) and $\log \|x\|$ (right) of the system \eqref{system_exp-1}. The parameter used are $t_0=0$, $x_0=(0.5,1)$ and $\alpha = 0.5$. \label{fig3}}
\end{center}
\end{figure}
\begin{figure}
\begin{center}
\includegraphics[scale=0.4]{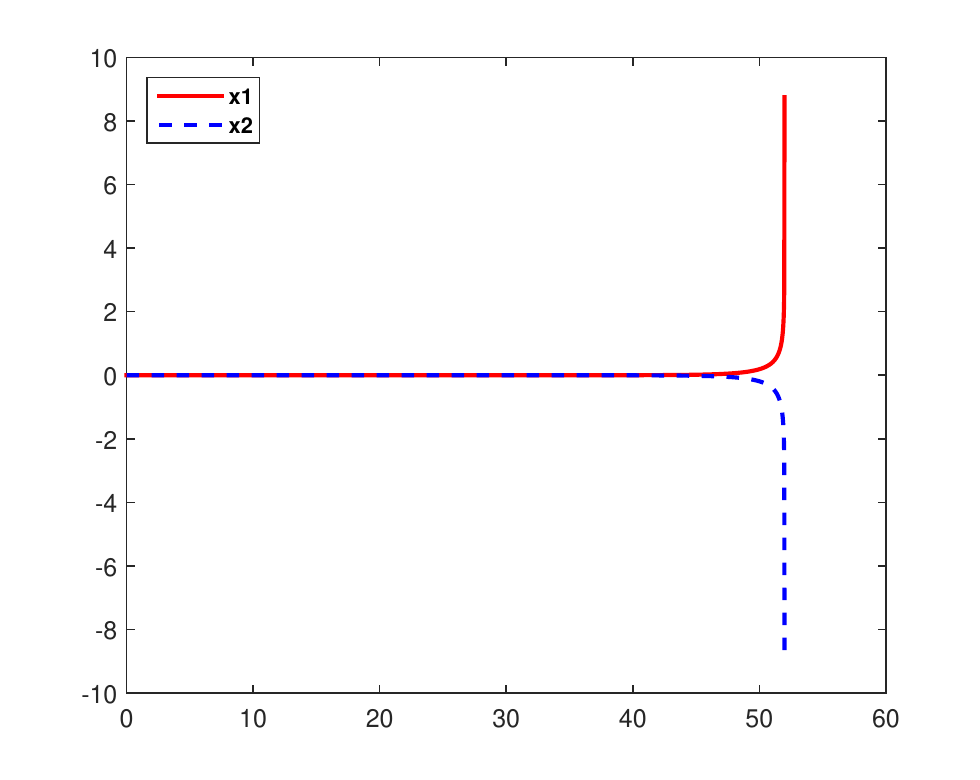}
\caption{Time evolution of the states $(x_1,x_2)$ of the system \eqref{system_exp-1}. The parameter used are $t_0=0$, $x_0=(10^{-10},10^{-9})$ and $\alpha = 1$. \label{fig4}}
\end{center}
\end{figure}
\end{ex}
%
%
\begin{ex}Consider the following nonlinear system
\begin{align}\label{system_sin5} 
\begin{cases}
\dot{x}_1=-x_1 + \frac{t\,\sqrt{\|x\|}}{1+t^2} + \psi(t)\\
\dot{x}_2=-x_2 + \frac{\sqrt{\|x\|}}{1+t^2} + \psi(t)
\end{cases}
\end{align}
with $\psi(t) = \frac{1}{\log(2+t)}$. The system \eqref{system_sin5} can be written as
\begin{equation}\label{system_sin_mat8}
\dot{x}=Ax+g(t,x),
\end{equation}
where $x=(x_1,x_2)\in \mathbb{R}^2$ is the state of the system, $A=-I_2$ and
\begin{equation*}
g(t,x) = \frac{\sqrt{\|x\|}}{1+t^2}(t,\,1) + \psi(t)(1,1).
\end{equation*}
The nominal system is globally uniformly exponentially stable. Moreover, we have
\begin{align*}
\|g(t,x)\|^2 & = \left(\frac{t\,\sqrt{\|x\|}}{1+t^2} + \frac{1}{\log(2+t)}\right)^2 + \left(\frac{\sqrt{\|x\|}}{1+t^2} + \frac{1}{\log(2+t)}\right)^2 \\
& \leq 2\left(\frac{t^2\,\|x\|}{(1+t^2)^2} + \frac{1}{\log^2(2+t)} + \frac{\|x\|}{(1+t^2)^2} + \frac{1}{\log^2(2+t)}\right) \\
& \leq 2\left(\frac{\|x\|}{1+t^2} + \frac{2}{\log^2(2+t)}\right).
\end{align*}
Hence $\|g(t,x)\| \leq \phi(t)\omega(\|x\|)+\lambda(t)$ with $\phi(t) = \frac{\sqrt{2}}{\sqrt{1+t^2}}$, $\omega(x) = \sqrt{\|x\|}$ and $\lambda(t) = \frac{2}{\log(2+t)}$. All the assumptions of Theorem \ref{theorem2_-1} (or Corollary \ref{cor_2}) are satisfied, then the system \eqref{system_sin5}  is globally uniformly practically exponentially stable. The figure \ref{fig5} shows the time evolution of the states $(x_1,\,x_2)$ of the system \eqref{system_exp-1} with the initial states $(x_{1,0},\,x_{2,0}) = (4,-1)$. One can notice that the solutions are stable but the origin is not attractive. This is in total agreement with Theorem \ref{theorem2_-1} since the function $\lambda$ do not belong to $L^p(\R_+)$ for any $p\in [1,\infty)$ but $\lim_{t\to \infty} \lambda(t)=0$. (see Theorem \ref{theorem2_-1} (iii)). Moreover, if we denote $r = \lim_{t\to \infty}\|x(t)\|$ then $\log(\|x(t)\|-r)$ decreases linearly as $t$ increases, which confirms the exponential decay of the solution toward $r$. Finally, if we consider the system \eqref{system_sin5} with $\psi(t) = \frac{1}{2+t}$ then using a same reasoning as here above we expect the origin to be attractive (since $t\mapsto\frac{1}{2+t}\in L^2(\R_+)$) and the system \eqref{system_sin5} is globally uniformly practically exponentially stable.  The figure \ref{fig6} confirms our study and shows that the solution goes to zero as $t$ tends to infinity.
\begin{figure}
\begin{center}
\includegraphics[scale=0.3]{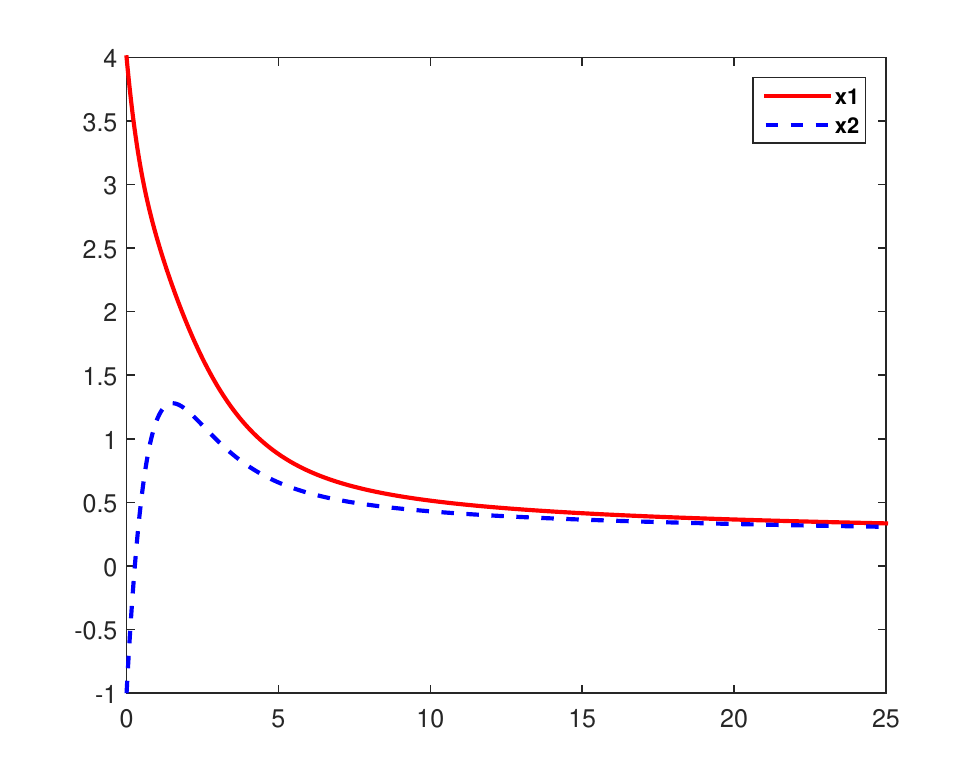} \includegraphics[scale=0.3]{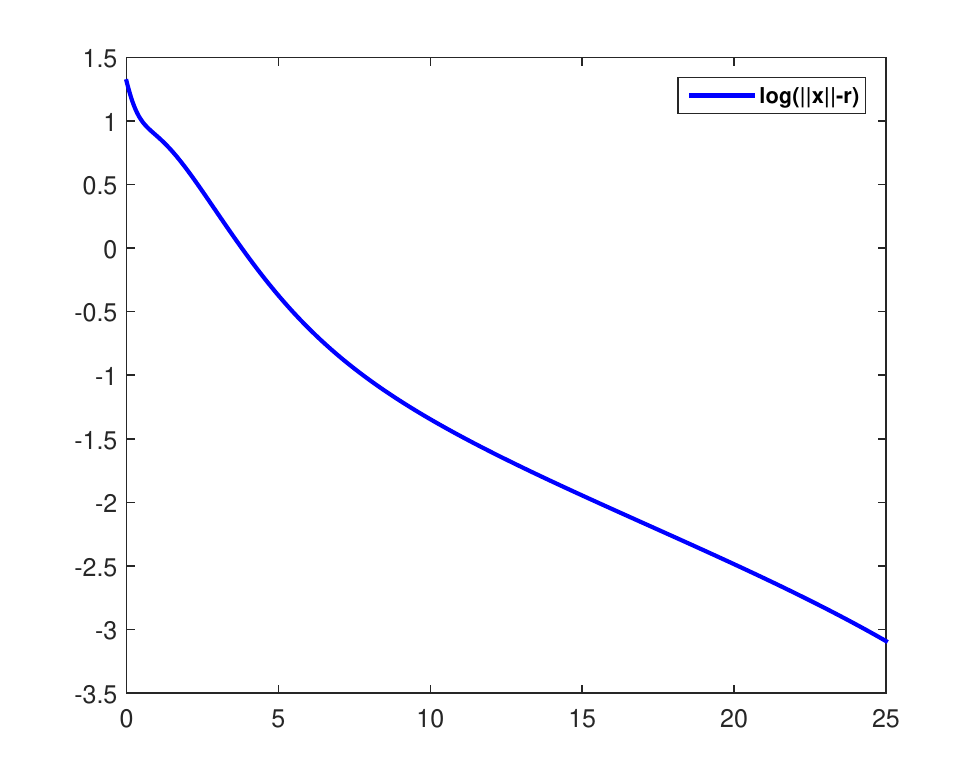}
\end{center}
\caption{Time evolution of the states $(x_1,x_2)$ (left) and $\log(\|x(t)\|-r)$ (right) of the system \eqref{system_sin5} with $\psi(t) = \frac{1}{\log(2+t)}$. The parameter used are $t_0=0$ and $x_0=(4,-1)$. \label{fig5}}
\end{figure}
\begin{figure}
\begin{center}
\includegraphics[scale=0.4]{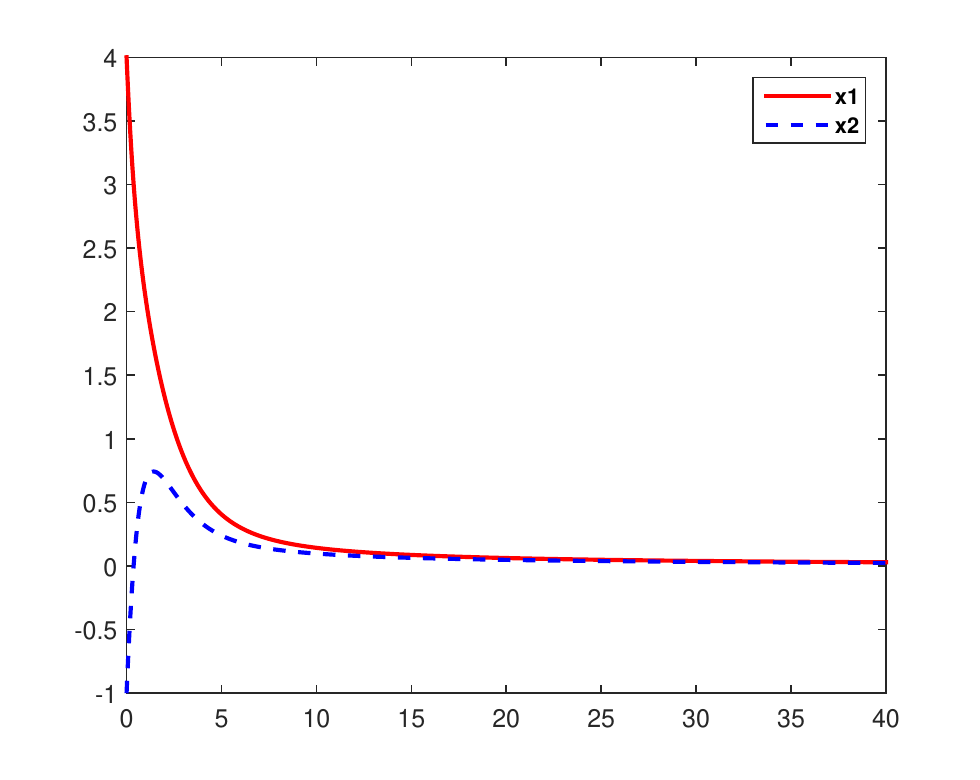}
\end{center}
\caption{Time evolution of the states $(x_1,x_2)$ of the system \eqref{system_sin5} with $\psi(t) = \frac{1}{2+t}$. The parameter used are $t_0=0$ and $x_0=(4,-1)$. \label{fig6}}
\end{figure}
\end{ex}

\section*{Conclusion}
We derived some sufficient conditions that ensure the global uniform exponential stability, the uniform exponential stability and the global uniform practical exponential stability of a perturbed system when the behavior of the perturbation terms is known (or can be estimated). Our approach is based on the Gronwall-Bellman inequalities instead of the Lyapunov techniques, which makes it easy to apply in practice. Several examples, especially when the standard Lyapunov approaches may fail, are given in illustration.



\end{document}